\documentclass[11pt]{article}
\usepackage{epsf}
\usepackage{color}
\usepackage{epsfig}

\usepackage{graphicx}
\usepackage{amsmath}
\usepackage{latexsym}
\usepackage{amssymb}
\usepackage{amsfonts}
\usepackage{hyperref}

\usepackage{a4}

\newcommand{\e}{{\rm e}}

\newcommand{\psis}{{\mathcal S}}

\newcommand{\psib}{\chi}

\begin{document}

\title{Splitting methods with complex coefficients}

\author{Sergio Blanes$^{1}$\thanks{Email: \texttt{serblaza@imm.upv.es}}
   \and
  Fernando Casas$^{2}$\thanks{Email: \texttt{Fernando.Casas@uji.es}}
   \and
 Ander Murua$^{3}$\thanks{Email: \texttt{Ander.Murua@ehu.es}}
   }

\maketitle


\begin{abstract}

Splitting methods for the numerical integration of differential
equations of order greater than two involve necessarily negative
coefficients. This order barrier can be overcome by considering
complex coefficients with positive real part. In this work we
review the composition technique used to construct methods of this
class, propose new sixth-order integrators and analyze their main
features on a pair of numerical examples, in particular how the
errors are propagated along the evolution.

\vspace*{0.5cm}

\begin{description}
 \item $^1$Instituto de Matem\'atica Multidisciplinar,
  Universidad Polit\'ecnica de Valencia, E-46022  Valencia, Spain.
 \item $^2$Institut de Matem\`atiques i Aplicacions de Castell\'o and Departament de Matem\`atiques,
 Universitat Jaume I,
  E-12071 Castell\'on, Spain.
 \item $^3$Konputazio Zientziak eta A.A. saila, Informatika
Fakultatea, EHU/UPV, Donostia/San Sebasti\'an, Spain.
\end{description}

\end{abstract}

\section{Introduction}

Splitting methods for the numerical integration of differential
equations constitute an appropriate choice when the associated vector field
can be decomposed into several pieces and
each of them is explicitly integrable.

Given the initial value problem
\begin{equation}   \label{eq.1.1}
   x' = f(x), \qquad x_0 = x(0) \in \mathbb{R}^D
\end{equation}
with $f: \mathbb{R}^D \longrightarrow  \mathbb{R}^D$ and solution
$\varphi_t(x_0)$, let us suppose that  $f$  can be expressed as $f
= \sum_{i=1}^m f^{[i]}$ for certain functions $f^{[i]}:
\mathbb{R}^D \longrightarrow  \mathbb{R}^D$, in such a way that
the equations
\begin{equation}   \label{eq.1.2}
   x' = f^{[i]}(x), \qquad x_0 = x(0) \in \mathbb{R}^D, \qquad i=1, \ldots, m
\end{equation}
can be integrated exactly, with solutions $x(h) =
\varphi_h^{[i]}(x_0)$ at $t = h$, the time step. Splitting methods
intend to approximate the exact flow $\varphi_h$ by a composition of flows
$\varphi_h^{[i]}$.
For instance,
\begin{equation}   \label{eq.1.3}
   \psib_h = \varphi_h^{[m]} \circ \cdots \circ \varphi_h^{[2]} \circ
   \varphi_h^{[1]},
   \qquad
   \psib_h^* = \varphi_h^{[1]} \circ \varphi_h^{[2]} \circ  \cdots \circ \varphi_h^{[m]}
\end{equation}
both provide first-order approximations to the exact
solution, since $\psib_h(x_0) = \varphi_h(x_0) + \mathcal{O}(h^2)$ and
similarly for $\psib_h^*$ (which is called the adjoint of $\psib_h$ and verifies
$\psib_h^*=\psib_{-h}^{-1}$).

It is possible to get higher order
approximations by introducing more maps with additional real
coefficients, $\varphi_{a_{ij} h}^{[i]}$,  in  (\ref{eq.1.3}). Perhaps the most popular splitting method is the second order symmetric composition
\begin{equation}  \label{eq:S_h}
\psis_h^{[2]} = \chi_{h/2} \circ \chi_{h/2}^* =
 \varphi_{h/2}^{[m]} \circ \cdots \circ \varphi_{h/2}^{[2]} \circ  \varphi_h^{[1]}
 \circ \varphi_{h/2}^{[2]} \circ  \cdots \circ \varphi_{h/2}^{[m]}.
\end{equation}
 When $f$ in (\ref{eq.1.1}) is separable in two parts the above particularizes to
\begin{equation}   \label{eq.1.4}
   \psib_h = \varphi_h^{[2]} \circ \varphi_h^{[1]}, \qquad
   \psib_h^* = \varphi_h^{[1]} \circ \varphi_h^{[2]}, \qquad
\psis_h^{[2]} =\varphi_{h/2}^{[2]} \circ  \varphi_h^{[1]} \circ
\varphi_{h/2}^{[2]},
\end{equation}
and $\psis_h^{[2]}$  is known as the
Strang splitting \cite{strang68otc}, the leapfrog or the
St\"ormer--Verlet method \cite{verlet67ceo}, depending on the
context where it is used. More generally, one may choose the coefficients $a_i$, $b_i$
to achieve order $r$ with the composition
\begin{equation}   \label{eq.1.AB}
\psi_h  =  \varphi^{[2]}_{b_{s+1} h}\circ \varphi^{[1]}_{a_{s}h}
   \circ \varphi^{[2]}_{b_{s} h}\circ \cdots \circ
 \varphi^{[2]}_{b_{2}h}\circ \varphi^{[1]}_{a_{1}h} \circ
 \varphi^{[2]}_{b_{1}h}.
\end{equation}
It turns out that $\psi_h$ can also be written in terms
of $\chi_{h}$ and $\chi_{h}^*$
\begin{eqnarray}
\psi_h & = &
  \left(\varphi^{[2]}_{\alpha_{2s}h} \circ \varphi^{[1]}_{\alpha_{2s}h}\right)
  \circ \cdots \circ
  \left(\varphi^{[2]}_{\alpha_{2}h} \circ \varphi^{[1]}_{\alpha_{2}h}\right)
  \circ
 \left(\varphi^{[1]}_{\alpha_{1}h} \circ \varphi^{[2]}_{\alpha_{1}h}\right)
  \nonumber   \\
 & = &  \chi_{\alpha_{2s} h}\circ
 \chi^*_{\alpha_{2s-1}h}\circ \cdots\circ \chi_{\alpha_{2}h}\circ
 \chi^*_{\alpha_{1}h}   \label{eq.1.MetAdj}
\end{eqnarray}
as long as
\begin{equation} \label{eq.1.RelCoefs}
a_j=\alpha_{2j-1}+\alpha_{2j}, \qquad  b_{j+1} = \alpha_{2j} +
\alpha_{2j+1}.
\end{equation}
Equivalently,
\begin{equation}\label{alfas_a_b}
 \alpha_{1}= b_1, \qquad  \alpha_{2j+1}= b_1 +\sum_{k=1}^{j}(b_{k+1}-a_k), \qquad
 \alpha_{2j}= \sum_{k=1}^{j}(a_k-b_{k}),
\end{equation}
with $\alpha_0=\alpha_{2s+1}=0$. This relation remains valid if $\sum_{i=1}^{s} a_i = \sum_{i=1}^{s+1} b_i$
 \cite{mclachlan95otn}. A relevant consequence of this property is that, starting with the
 coefficients $a_i,b_i$ of a given splitting method,  we can get the
coefficients $\alpha_i$ for the composition (\ref{eq.1.MetAdj}), which can be then applied
in a more general setting with the maps
$\chi_{h}$ and $\chi_{h}^*$ of (\ref{eq.1.3}). A particular case widely used in practice to achieve
high order approximations consists in considering compositions using the Strang splitting
(\ref{eq:S_h}) as basic method,
\begin{equation}\label{CompS2}
\psi_h =
 \psis^{[2]}_{\alpha_{s}h}\circ \cdots\circ \psis^{[2]}_{\alpha_{2}h}\circ
 \psis^{[2]}_{\alpha_{1}h}.
\end{equation}

Splitting methods are, in general, explicit, easy to implement and
preserve structural properties of the exact solution, thus
conferring to the numerical scheme a qualitative superiority with
respect to other standard integrators, especially when long time
intervals are considered (see \cite{blanes08sac} for a review).
Examples of these structural features are symplecticity, volume
preservation, time-symmetry and conservation of first integrals.
In this sense, splitting methods constitute an important class of
\emph{geometric numerical integrators}
\cite{hairer06gni,iserles00lgm,leimkuhler04shd,mclachlan02sm,mclachlan06gif,sanz-serna94nhp}.

It has been shown that some of the coefficients in splitting schemes (\ref{eq.1.AB})
are negative when the order $r \ge 3$ \cite{goldman96noo,sheng89slp,suzuki91gto}.
In other words, the methods always involve stepping backwards in time.
An elementary proof of this feature can be worked out as follows. It is quite straightforward to
check that one of the necessary condition for the composition
(\ref{eq.1.MetAdj}) (respectively, (\ref{CompS2})) to have order $r \ge 3$ is
\begin{equation}\label{Cond3}
\sum_{i=1}^{k} \alpha_i^3 = 0,
\end{equation}
with $k=s$ (respect. $k=2s$). Obviously, this sum vanishes only  if at least
one of the  $\alpha_i$ is negative. In consequence, the flows
$\varphi_h^{[j]}, \ j=1,\ldots,m-1$ in (\ref{eq:S_h}) and
$\varphi_h^{[j]}, \ j=2,\ldots,m-1$ in  (\ref{eq.1.3}) evolve with
at least one negative fractional time step. On the other hand, by taking into account
the link (\ref{eq.1.RelCoefs}) among coefficients of  (\ref{eq.1.AB}) and (\ref{eq.1.MetAdj}),
condition (\ref{Cond3}) with $k=2s$ leads to
\[
 \sum_{j=1}^{s}(\alpha_{2j-1}^3+\alpha_{2j}^3)=0 \quad \Rightarrow
  \quad  \exists \ k \ / \ \alpha_{2k-1}^3+\alpha_{2k}^3 < 0 \quad \Rightarrow
  \quad    a_k=\alpha_{2k-1}+\alpha_{2k} < 0.
\]
In a similar way, using the same condition with
$\alpha_0=\alpha_{2s+1}=0$, one has
\[
 \sum_{j=0}^{s} (\alpha_{2j}^3 + \alpha_{2j+1}^3)=0
  \quad  \Rightarrow
  \quad   \exists \ l \ / \ \alpha_{2l}^3 + \alpha_{2l+1}^3 < 0
  \quad \Rightarrow  \quad    b_l=\alpha_{2l} + \alpha_{2l+1} < 0
\]
and then at least one $a_i$ as well as one $b_i$ are negative.
It must be stressed that condition (\ref{Cond3}) still persists when
the processing technique is used, so that the same conclusion also follows  in this case
\cite{blanes05otn}.

In summary, the presence of negative coefficients in splitting
methods of order higher than two is unavoidable if one restricts oneself to
\emph{real} coefficients. Of course, this does not suppose any special impediment
when the flow of the ODE
evolves in a group (such as in the
Hamiltonian case), but may be unacceptable when the differential equation is
defined in a semigroup \cite{mclachlan02sm}, as occurs, for instance, with
the simple heat equation $u_t = \Delta u$ on the unit interval with
homogeneous Dirichlet boundary conditions. Then the corresponding generated semigroup
is well defined only for $t \ge 0$ \cite{hansen09hos}.

More generally, consider the nonlinear heat equation
 \begin{equation}\label{NL_heat2}
  \frac{\partial}{\partial t} u(x,t) =
  \sum_{i=1}^dD_i(v_i(x)D_i u(x,t)) + F(x,u(x,t))
 \end{equation}
with functions $v_i$ real and positive, and $D_i \equiv \partial /\partial x_i$, on a certain domain $\Omega \in \mathbb{R}^d$. If a space discretization is carried out (either
by finite differences or by a pseudospectral method), a large system of ODEs results
which has to be numerically integrated in time. To this end, we can split the resulting equation into
linear and nonlinear parts, but
schemes of the form
(\ref{eq.1.MetAdj}) or (\ref{CompS2}) of at most order $r=2$ can only be applied,
since the resulting discrete Laplacian with negative fractional time steps is not well conditioned.

A closely related problem is
the linear Schr\"odinger equation ($\hbar = 1$):
\begin{equation}\label{Schrod1}
  i\frac{\partial}{\partial t} \Psi(x,t) =
  \left( -\frac{1}{2m} \Delta + V(x) \right)  \Psi(x,t).
\end{equation}
A technique used in practice to obtain the eigenvalues and eigenfunctions for a given potential $V$
consists in numerically integrating the equation (after spatial discretization)
along pure imaginary times ($\tau=-it$). Equivalently,  the equation to be analyzed is
\begin{equation}\label{Schrod2}
  \frac{\partial}{\partial \tau} \Psi =
  \left( \frac{1}{2m} \Delta - V(x) \right)  \Psi,
\end{equation}
which can be considered as a linear heat equation. The system
evolves to the ground state whose norm decreases exponentially in
proportion to the value of its energy (eigenvalue). By
orthogonalization, one can make the system to evolve to any other
eigenfunction \cite{auer01afo,lehtovaara07sot}. In any case,
whereas there is no special difficulties with numerically
integrating equation (\ref{Schrod1}) using a splitting methods
with negative fractional time steps, this is not the case for
(\ref{NL_heat2}) and (\ref{Schrod2}) due to the presence of the
Laplacian.

It has been noticed, however, that  higher order splitting methods
with complex coefficients having positive real part do exist
\cite{bandrauk91ies,mclachlan02sm,suzuki90fdo,suzuki91gto,suzuki95hep}.
These schemes were reported mainly for theoretical purposes but
received very little attention as practical numerical tools. Perhaps the main reason was
that working with complex arithmetic makes
the schemes more involved and, in many cases, also considerably more
costly from a computational point of view (usually, four times more expensive).

It is only within recent years that a systematic
search for new methods  with complex coefficients has been carried out and the resulting
schemes have been tested in different settings:
Hamiltonian systems in celestial
mechanics  \cite{chambers03siw}, the time-dependent
Schr\"odinger equation in quantum mechanics
\cite{bandrauk06cis,prosen06hon} and also in the more abstract setting of
evolution equations with unbounded operators generating analytic
semigroups \cite{castella09smw,hansen09hos}. In this sense, we recall that the
propagator $\exp(z \Delta)$ ($z \in \mathbb{C}$) associated with the Laplacian
 is well defined (in a reasonable distributional sense) if and only if $\mbox{Re}(z) \ge 0$
 \cite{castella09smw}. More generally, it is possible to extend the semigroup related with
 parabolic PDEs into a sector in the right half plane of $\mathbb{C}$ \cite{hansen09hos}.


The aim of this paper is to review some of the splitting methods with complex coefficients
published in the literature, propose new sixth-order schemes in the class (\ref{CompS2}) and
analyze them on a pair of simple numerical examples, to get a glance of the performance
and main features of this kind of integrators and some of the difficulties involved.

\section{Integrators with complex coefficients} \label{ntsteps}


Most of the existing splitting methods with complex coefficients  have been
constructed by applying the composition technique to the symmetric second-order
leapfrog scheme $\psis^{[2]}$. Thus, one gets a third-order method as
\begin{equation}  \label{eq:S_h-C3}
   \psis_{h}^{[3]} =
   \psis_{\alpha h}^{[2]} \circ \psis_{\beta h}^{[2]},
\end{equation}
where the coefficients have to satisfy (\ref{Cond3}) together with the consistency condition
\[
\left. \begin{array} {rcl}
   \alpha + \beta & = & 1 \\
   \alpha^3 + \beta^3 & = & 0
\end{array} \right\}
\Rightarrow  \;\;
\begin{array}{l}
  \alpha = \frac{1}{2} \mp i\frac{\sqrt{3}}{6}, \qquad
   \beta = \frac{1}{2} \pm i\frac{\sqrt{3}}{6}.
\end{array}
\]
Due to its simplicity, this scheme has been rediscovered several
times, either as the composition (\ref{eq:S_h-C3})
\cite{bandrauk91ies,suzuki90fdo,castella09smw} or by solving the order
conditions required by (\ref{eq.1.AB}) with $s=2$
\cite{chambers03siw,hansen09hos}.

A fourth-order integrator can be obtained with the symmetric
composition
\begin{equation}\label{suzu1}
   \psis_{h}^{[4]} =
   \psis_{\alpha h}^{[2]} \circ \psis_{\beta h}^{[2]} \circ  \psis_{\alpha
   h}^{[2]}.
\end{equation}
Although the necessary order conditions are the same, the time-symmetry of the composition rises
the order by one (all the error terms at odd orders vanish identically):
\[
\left. \begin{array} {rcl}
   2\alpha + \beta & = & 1 \\
   2\alpha^3 + \beta^3 & = & 0
\end{array} \right\}
\Rightarrow  \;\; \alpha = \frac{1}{2-2^{1/3}\, \e^{2ik\pi/3}}, \qquad
 \beta = \frac{2^{1/3}\, \e^{2ik\pi/3}}{2-2^{1/3}\, \e^{2ik\pi/3}}
\]
with $k=0,1,2$. Notice that for $k=1,2$ it is true that
$\mbox{Re}(\alpha), \mbox{Re}(\beta) >0$.

Another fourth-order method can be obtained by symmetrizing
the third-order scheme (\ref{eq:S_h-C3}), i.e.,
\begin{equation}  \label{eq:S_h-C34}
   \psis_{h}^{[4]} =
   \psis_{\alpha/2 h}^{[2]} \circ \psis_{\beta/2 h}^{[2]} \circ
   \psis_{\beta/2 h}^{[2]} \circ \psis_{\alpha/2 h}^{[2]}.
\end{equation}

Methods (\ref{eq:S_h-C3}), (\ref{suzu1}) and (\ref{eq:S_h-C34}) can be used to generate
recursively higher
order composition schemes as
\begin{equation}  \label{eq:S_h-C3n}
   \psis_{h}^{[n+1]} =
   \psis_{\alpha h}^{[n]} \circ \psis_{\beta h}^{[n]}.
\end{equation}
Here the coefficients have to verify the conditions $\alpha + \beta = 1$, $\alpha^{n+1} + \beta^{n+1}  =  0$,
whence
\[
\begin{array}{l}
  \alpha = \displaystyle \frac{1}{2} +
  i \, \frac{\sin (\frac{2l+1}{n+1} \pi)}{2+2\cos (\frac{2l+1}{n+1} \pi )} \quad
   \mbox{for} \quad
   \left\{ \begin{array} {ll}
   -\frac{n}{2} \leq l \leq \frac{n}{2}-1 & \mbox{if} \ n \ \mbox{is even}, \\
   -\frac{n+1}{2} \leq l \leq \frac{n-1}{2} & \mbox{if} \ n \ \mbox{is odd},
    \end{array} \right.
\end{array}
\]
and $\beta=1-\alpha$. The choice $l=0$ gives the solutions with
the smallest phase and allows one to build methods up to order six with coefficients having positive real part.
This feature was stated in
\cite{suzuki95hep} and rediscovered in \cite{hansen09hos}.

In a similar way, one may use recursively a symmetric three term composition,
which allows to increase the order by two at each iteration:
\begin{equation}\label{suzu_n}
   \psis_{h}^{[n+2]} =
   \psis_{\alpha h}^{[n]} \circ \psis_{\beta h}^{[n]} \circ  \psis_{\alpha
   h}^{[n]},
\end{equation}
with
\[
   2\alpha + \beta  =  1, \qquad
   2\alpha^{n+1} + \beta^{n+1}  =  0.
\]
The solutions providing coefficients with the smallest
phase are
\[
 \alpha =\frac{\e^{i\pi/(n+1)}}{2^{1/(n+1)}-2 \, \e^{i\pi/(n+1)}}, \qquad
 \beta = 1-2\alpha,
\]
and methods up to order eight with coefficients having positive
real part are possible. Moreover, methods up to order fourteen of
the more general form (\ref{CompS2}) with coefficients $\alpha_j$
with positive real part are attainable
\cite{castella09smw,hansen09hos}. An interesting (and open)
question is to determine whether arbitrarily high orders can be
attained or wether, as for the previous compositions, there is an
order barrier for methods of the form (\ref{CompS2}) with
$\mathrm{Re}(\alpha_j)>0$. Observe that any method of the form
(\ref{CompS2}) with coefficients having positive real part can be
expressed in terms of the elementary flows $\varphi_{\lambda
h}^{[j]}$ with $\mathrm{Re}(\lambda)>0$ when $\mathcal{S}_h^{[2]}$
is taken as the leapfrog  (\ref{eq.1.4}). This is also true, of
course, in the more general case when $f$ is split in $m$ parts
and  $\mathcal{S}_h^{[2]}$ is  taken as the symmetric second order
basic method (\ref{eq:S_h}).

For instance, suppose that $f$ in (\ref{eq.1.1}) is separable in two parts, so that
$\mathcal{S}_h^{[2]}$ is given by (\ref{eq.1.4}). Then it is straightforward to check that the third order
scheme (\ref{eq:S_h-C3}) can be written as
\begin{equation}    \label{chambers.1}
     \mathcal{S}_h^{[3]} = \varphi_{b_3 h}^{[2]} \circ \varphi_{a_2 h}^{[1]} \circ \varphi_{b_2 h}^{[2]}
      \circ \varphi_{a_1 h}^{[1]} \circ \varphi_{b_1 h}^{[2]}
\end{equation}
with $a_1=  \frac{1}{2} + i\frac{\sqrt{3}}{6}, \ a_2=a_1^*, \ b_1
= a_1/2, \ b_2=1/2, \ b_3=b_1^*$.  This particular symmetry of the coefficients results in a method
whose leading error terms at order 4 are all strictly imaginary \cite{chambers03siw}.

Another question of practical nature is the construction of methods of the form (\ref{CompS2})
with $\mathrm{Re}(\alpha_j)>0$ involving the minimum number $s$ of compositions for a prescribed order. For instance, the minimal number of compositions for achieving order 6 is $s=7$. The corresponding
order conditions can be written as \cite{blanes08sac,hairer06gni,murua99ocf}
\begin{eqnarray}
  \label{eq:order6}
&&  \sum_{j=1}^{s} \alpha_j=1,\quad
  \sum_{j=1}^{s} \alpha_j^{k}=0,\quad k=3,5, \\
&& \sum_{j=1}^{s} \alpha_j^{k}c_j^{\ell}=0,\quad (k,\ell)\in \{(3,1), (3,2),(3,3),(5,1)\},
\end{eqnarray}
where for each $j=1,\ldots,s$,
\begin{eqnarray*}
  c_j = \frac{\alpha_j}{2} + \sum_{i=1}^{j-1} \alpha_i.
\end{eqnarray*}
This system of algebraic equations  has several solutions with
$\mathrm{Re}(\alpha_j)>0$. Among them, we have chosen the two sets
of coefficients collected in Table~\ref{table.coefs1}.  The first
one corresponds to a symmetric method, $\alpha_{s+1-i} =
\alpha_i$, as scheme (\ref{suzu1}), and was already found by
Chambers \cite{chambers03siw}. The second method is apparently
new, and possesses the special symmetry $\alpha_{s+1-i} =
\alpha_i^*$, as scheme (\ref{eq:S_h-C3}) (or (\ref{chambers.1})
when expressed as (\ref{eq.1.AB})).

\begin{table}[htb]
\begin{center}
\caption{Coefficients of two 7-stage sixth-order methods of type (\ref{CompS2}): S$_{7}6$
is a symmetric method and S$_{7}^*6$ is conjugate to a symmetric
method (symmetric in the real part of the coefficients and
skew-symmetric in the imaginary part).}
\label{table.coefs1}
\vspace*{0.3cm}
\begin{tabular}{l}
\begin{tabular}{l}
\hline
  S$_{7}6$   \\
\hline
 $\alpha_1= 0.116900037554661284389 + 0.043428254616060341762 \,i$  \\
 $\alpha_2= 0.12955910128208826275 - 0.12398961218809259330 \,i$  \\
 $\alpha_3= 0.18653249281213381780 + 0.00310743071007267534 \,i$ \\
 $\alpha_{4}=  0.13401673670223327014 + 0.15490785372391915239 \,i$ \\
 $\alpha_5= \alpha_3, \quad \alpha_6= \alpha_2, \quad \alpha_7= \alpha_1  $  \\
\hline \hline
   S$_{7}^*6$    \\
\hline
 $\alpha_1=  0.133741778914683628452 - 0.028839028371025553995 \, i$ \\
 $\alpha_2= 0.12134019583938803504 + 0.11585180844272788007 \, i$  \\
 $\alpha_3= 0.13489797942731665044 -  0.12906241362827633477 \, i$ \\
 $\alpha_{4}= 0.22004009163722337213$ \\
 $\alpha_5= \alpha_3^*, \quad \alpha_6= \alpha_2^*, \quad \alpha_7= \alpha_1^*  $  \\
\hline
\end{tabular}
\end{tabular}
\\
\end{center}
\end{table}

\section{Numerical examples}

\subsection{Example 1: the harmonic oscillator}

We consider the simple harmonic oscillator to illustrate some qualitative properties of
the previous composition methods with complex coefficients. That is, we take
the Hamiltonian
function $H(q,p) = \frac{1}{2}(p^2 + q^2)$, with $q,p \in
\mathbb{R}$. The corresponding equations of motion are
linear  and can be written as
\begin{equation}\label{harmonic2}
    x' \equiv \left(\begin{array}{c}
        q^\prime \\
        p^\prime  \end{array}\right) =
    \Big[   \underbrace{\left(\begin{array}{cc}
        0 & 1 \\
        0 & 0  \end{array}\right)}_{A} +
         \underbrace{\left(\begin{array}{cc}
        0 & 0 \\
       -1 & 0  \end{array}\right)}_{B}
 \Big]
     \left(\begin{array}{c}
        q \\
        p  \end{array}\right) = (A + B) \, x,
\end{equation}
so that the numerical solution at time $t=h$ furnished by method
(\ref{eq.1.AB})  is given by
\begin{eqnarray}
  \label{x(h)harmonic}
     x(h) = K(h) x_0 \equiv \e^{b_{s+1} h B} \, \e^{a_s h A} \, \e^{b_s h B} \, \cdots
\e^{b_2 h B} \, \e^{a_1 h A} \, \e^{b_{1} h B} x_0.
\end{eqnarray}
As is well known, for splitting methods with real coefficients the
average error in energy remain constant for exponentially long
times under suitable general conditions on the Hamiltonian. For
the particular case of the harmonic oscillator and with a
sufficiently small time step, this is true for all times,  and the
average error in positions grows only linearly.

\begin{figure}[h!]
\begin{center}
\makebox{\epsfig{figure=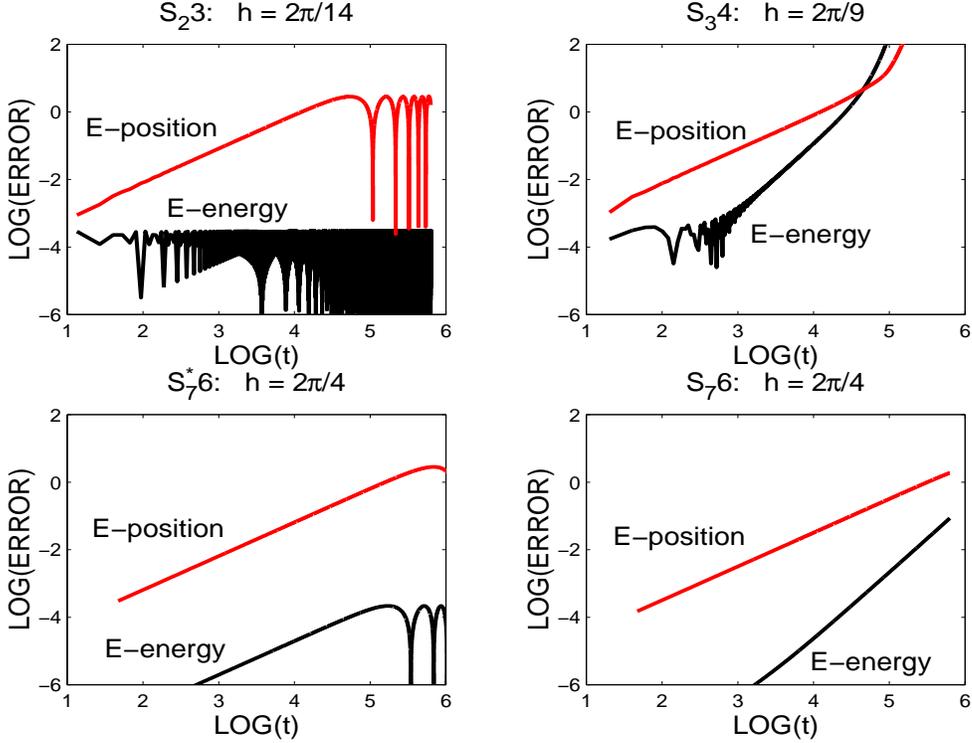,height=10cm,width=13cm}}
\end{center}
\caption{{Error in position and energy (taking the real part from
the output) obtained with the 4th- and 6th-order symmetric schemes $S_34$ and $S_76$, and
non-symmetric methods $S_23$ and $S_7^*6$.
The time step is chosen such
that all schemes require 27-28 evaluations per period.}}
 \label{fig1}
\end{figure}

We propose here to check whether this also holds for methods with complex coefficients. To do that,
we take as initial conditions
$(q,p)=(1,1)$ and integrate the system (\ref{harmonic2}) for $t\in[0,20000 \pi]$ using a
constant time step. We measure the error in position and
energy of the output obtained by propagating the solution with the splitting method and then computing the real parts of the results
$q_{out}= \mbox{Re}(q)$, $p_{out}= \mbox{Re}(p)$. Figure~\ref{fig1} shows the results obtained
with the following
methods: (i) $S_23$, the 2-stage third-order non-symmetric method
(\ref{eq:S_h-C3}), (ii) $S_34$, the 3-stage fourth-order symmetric
method (\ref{suzu1}), (iii) $S_7^*6$, the 7-stage sixth-order
non-symmetric method, (iv) $S_76$,
the 7-stage sixth-order symmetric method. The coefficients of these two 6th-order methods are
collected in
Table~\ref{table.coefs1}. The time step is chosen such that all
methods require 27-28 evaluations per period. Notice the
significant difference in the qualitative behavior of the numerical solution. Whereas the error
grows exponentially for the symmetric methods $S_34$ and $S_76$, this is not the case
for $S_23$ and $S_7^*6$, which show a performance analogous to standard splitting methods with
real coefficients: bounded energy error and linear growth  of error in positions.
Of course, such a behavior deserves
a theoretical explanation, which we pursue next.

The matrix $K(h)$ in (\ref{x(h)harmonic}) is given explicitly
by
\[
   K(h) = \left( \begin{array}{cr}
                          1  &  0 \\
                         -b_{s+1} h & 1
                    \end{array} \right) \,
                \left( \begin{array}{cr}
                          1  &  a_s h \\
                          0 & 1
                    \end{array} \right) \,  \cdots
                \left( \begin{array}{cr}
                          1  &  a_1 h \\
                          0 & 1
                    \end{array} \right) \,
               \left( \begin{array}{cr}
                          1  &  0 \\
                         -b_{1} h & 1
                    \end{array} \right).
\]
In this way, one gets
\[
    K(h) =      \left( \begin{array}{rr}
                         p(h)+d(h)  &  \  q(h)+e(h) \\
                          -q(h)+e(h) & \  p(h)-d(h)
                    \end{array} \right)
\]
where $p(h)$, $d(h)$ (respectively, $q(h)$, $e(h)$) are
even (resp. odd) polynomial functions having in general complex coefficients
and $\det K(h) = p(h)^2+q(h)^2-d(h)^2-e(h)^2=1$.

If the splitting method  (\ref{eq.1.AB}) is such that

\begin{equation}\label{adjoint-symmetry}
 a_{s-j+1}=a_j^*, \qquad b_{s-j+2} = b_{j}^*
\end{equation}
(as happens, in particular, when it comes from a composition of
the form (\ref{CompS2}) with
  $\alpha_{s-j+1} = \alpha_{j}^*$),
then $K(h)^{-1} = K(-h)^*$. More specifically,
\[
    \left( \begin{array}{cr}
                           p(h) - d(h)  &  -q(h)-e(h) \\
                          q(h)-e(h) &   p(h)+d(h)
                    \end{array} \right) =
      \left( \begin{array}{cr}
               p(h)^*+d(h)^*  &  -q(h)^*-e(h)^* \\
             q(h)^*-e(h)^* &   p(h)^*-d(h)^*
                    \end{array} \right).
\]
This implies that $p(y)$, $q(h)$, and $e(h)$ are real polynomials, whereas the coefficients of $d(y)$ are purely imaginary. Notice that this is precisely the case of methods $S_23$ and $S_7^*6$.

If, on the other hand, the splitting method is symmetric, i.e., it is of the form (\ref{eq.1.AB}) satisfying
\begin{eqnarray*}
 a_{s-j+1}=a_j, \quad b_{s-j+2} = b_{j}
\end{eqnarray*}
(as happens, in particular, when it comes from a composition of the form (\ref{CompS2}) with
  $\alpha_{s-j+1} = \alpha_{j}$),
then $K(h)^{-1} = K(-h)$. This clearly implies that  $d(h)\equiv 0$, but in general the polynomials
$p(h)$, $q(h)$, and $e(h)$ have complex coefficients. For instance, methods
(\ref{suzu1}) ($S_34$) and $S_76$ are such that $p(h)$ has non-real coefficients.

When a splitting method with matrix $K(h)$ is used to integrate the harmonic oscillator, it is
essential that $p(h) \in \mathbb{R}$. Otherwise $K(h)^n$ grows exponentially with the number $n$ of steps.
As a matter of fact, the eigenvalues of $K(h)$ are
$\lambda_{1}=\e^{i\phi(h)}$ and $\lambda_{2}=\e^{-i\phi(h)}$, where
\begin{eqnarray*}
  \phi(h) = \arccos(p(h)),
\end{eqnarray*}
and thus $\max(|\lambda_1|,|\lambda_2|)>1$ if $p(h)\not \in \mathbb{R}$
(and also if $p(h) \in \mathbb{R}$ and $|p(h)|>1$).
That is precisely the situation with methods $S_34$ and $S_76$, and thus they are useless when integrating harmonic oscillators or systems that can be considered as close perturbations of harmonic oscillators with the partition (\ref{harmonic2}).

From the previous comments, it is clear that instability will take place when integrating the harmonic oscillator
unless $-1\leq p(h) \leq 1$. In fact, the numerical solution can still be (weakly) unstable when $p(h)^2=1$ with $q(h)d(h)e(h)\neq 0$ \cite{blanes08otl}. Furthermore,
 it is shown in  \cite{blanes08otl} that, for stable numerical solutions (that is, either $-1 < p(h) <1$ or $p(h)^2=1$ with $q(h)=d(h)=e(h)=0$), one has
\begin{eqnarray*}
     K(h)^n= Q(h)^{-1}
\left(
 \begin{array}{rr}
    \cos (n \phi(h)) &  \sin(n \phi(h))  \\
-\sin(n \phi(h)) &  \cos(n \phi(h))
 \end{array} \right) Q(h),
\end{eqnarray*}
with a suitable $2\times 2$ matrix $Q(h)$ (typically close to the identity matrix). In consequence,
the numerical solution $x_n=(q_n,p_n)$ is such that $\tilde x_n := Q(h) x_n$ corresponds to the exact solution at $t_n=n h$ of a harmonic oscillator with frequency $\tilde \omega = 1/h \phi(h) \approx 1$. This
feature explains why schemes $S_23$ and $S_7^*6$, when applied to the harmonic oscillator
(\ref{harmonic2})  with $h=\pi/7$ and $h= \pi/2$ respectively, exhibit a linear error growth in positions and a bounded error in energy, since for such methods,  $p(h)=1-h^2/2+\cdots$ is real and satisfies $p(h)  \in (-1,1)$
for the values of $h$ considered in the numerical experiments.

\subsection{Example 2: The Volterra--Lotka problem}

Consider now the Volterra--Lotka problem
\begin{equation}\label{volterra-lotka}
 \dot{u} = u(v-2), \qquad  \dot{v} = v(1-u).
\end{equation}
This is a very simple nonlinear system which allows us to make a
preliminary study about the behavior and performance of
methods with complex coefficients in the transition process from a linear
to a nonlinear problem. In a neighborhood of  the steady state at
$(u^*,v^*)=(1,2)$ the system can be considered close to a harmonic
oscillator. The nonlinear contributions are manifest as we move
away from it. The system evolves along periodic trajectories
around the equilibrium point in the region $0<u,v$ determined by the
first integral $I(u,v)=\ln (u v^2) - (u +v)$.

The vector field $f(u,v)=(u(v-2),v(1-u))$ can be separated in two
solvable parts and this can be done in different ways. We consider the
following split: $f_A=(u(v-2),0)$ and $f_B=(0,v(1-u))$ (although the linear
and nonlinear separation can also be considered).

\begin{figure}[h!]
\begin{center}
\makebox{\epsfig{figure=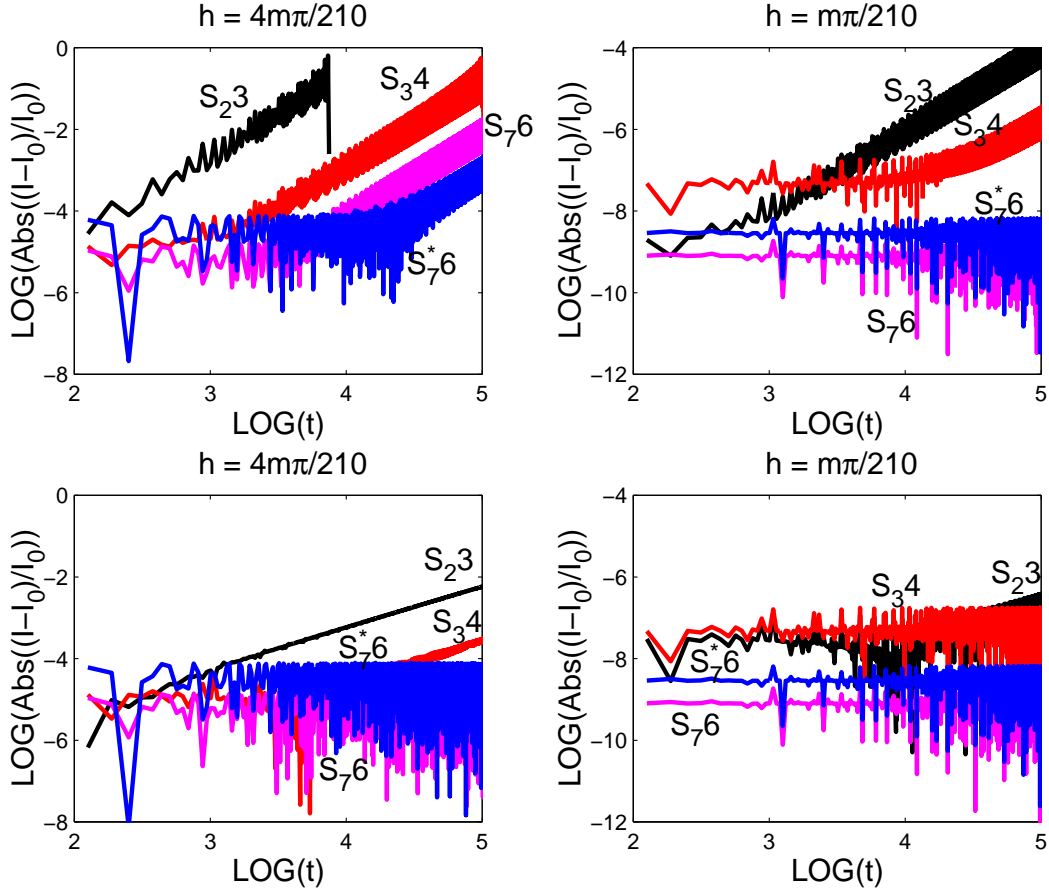,height=12cm,width=14cm}}
\end{center}
\caption{{Relative error in the first integral $I=\ln (u v^2) - (u
+v)$
 for the Volterra--Lotka problem with initial conditions $(u_0,v_0)=(2,4)$
 for the time steps  $h=\frac{4m\pi}{210}$ and
 $h=\frac{m\pi}{210}$, with $m$  the number of stages of each method.}}
 \label{fig3}
\end{figure}

We take as initial conditions $(u_0,v_0)=(2,4)$, integrate up to
$t=20000\times 2\pi$ and measure the relative error in the first
integral, $|I-I_0|/|I_0|$. As in the previous example, we
integrate using complex arithmetic and take the real parts of $u$
and $v$ only for representing the output. Figures~\ref{fig3}-(a)
and (b)  show the results obtained for time steps
$h=\frac{4m\pi}{210}$ and four times smaller $h=\frac{m\pi}{210}$,
with $m$ the number of stages of each method. In this way, all
methods require the same number of evaluations. Contrarily to the
pure harmonic oscillator, we observe a secular error growth in the
determination of the first integral for all methods which
diminishes considerably when the time step is reduced. The
observed behavior resembles what takes place with the so-called
pseudo-symplectic methods (integrators of order $n$ which preserve
symplecticity up to order $p>n$), where the dominant errors behave
as $\mathcal{E}_I=Ch^n+tDh^p$ for some constants $C$ and $D$. If
$p>n$ the secular part of the error does not manifest for
relatively long times when the time step is reduced.

We have repeated the same experiment but, after each time step, we discard the imaginary part
of $u$ and $v$ and initiate the next step only with their real part. In other words, we project each component
on the real axis at the end of each integration step. The results obtained are
shown in Figures~\ref{fig3}-(c) and (d).  Obviously, this way of proceeding does not preserve
symplecticity any more but the results obtained suggest that
a significant improvement in accuracy can be achieved.

\section{Conclusions and outlook}

We have presented a short review of the splitting and composition
technique to build methods of order greater than two with complex
coefficients with positive real part. This procedure allows to
overcome the order barrier where splitting methods of order
greater than two involve necessarily negative coefficients in the
real space. In general, splitting methods with complex
coefficients are considerably more expensive than the
corresponding methods with real coefficients (about four times more
expensive), and this make them hardly competitive in practice. For
this reason, one can think that the main application of the new
methods could be on parabolic PDEs, where higher order methods
with real coefficients (which necessarily have some negative coefficientes) can
not be used. However, there is a number of problems which evolve
in the complex space where using methods which complex
coefficients does not necessarily mean increasing the cost of
the algorithm. This can be the case, for instance, of the
Sch\"odinger equation (\ref{Schrod1}).

As for the practical implementation of splitting methods with
complex coefficients, in \cite{chambers03siw} it is claimed that
one has to carry the numerical integration in complex variables,
and (for problems with real solutions) one should take either the
real part of the variables or their modulus only for the output.
However, we have observed that removing the imaginary part at each
step, i.e. projecting on the real space at each step, the error
grow can be considerably diminished in some cases. In the
numerical examples considered in previous section, the linear
error grow in the first integrals originate from different sources
depending on wether the projection onto the real domain is
performed after each step or not. In the first case, the
projection after each step destroys symplecticity but only at a
higher order, and the schemes can be considered as
pseudosymplectic. In the second case, the method is actually
symplectic and can thus be (formally) considered as an exact
solution of a Hamiltonian system in the complex domain, which have
qualitatively different properties to trajectories in the real
domain. We have also noticed that the higher order methods present
a considerably reduced error grow. Then, it seems appropriate to
look for efficient higher order methods with complex coefficients.
In general, symmetric splitting methods are desirable. However, we
have shown that for the harmonic oscillator symmetric methods
(with non-real stability polynomial) present an exponential error
grow, which is not the case for methods with the special symmetry
(\ref{adjoint-symmetry}).
In a preliminary
search of methods, we have presented a new sixth-order method with
that special symmetry. This is an interesting subject to be
further explored since many problems in different applications can
be considered as perturbations to the harmonic oscillator.

\subsection*{Acknowledgements}

This work has been supported by Ministerio de Ciencia e Innovaci\'on
(Spain) under project MTM2007-61572 (co-financed by the ERDF of
the European Union).
SB also acknowledges financial support from Generalitat Valenciana
through project GV/2009/032.

\end{document}